\documentclass[12pt,leqno]{amsart}
\usepackage{amssymb} 
\begin{document}
\theoremstyle{plain}
\newtheorem{thm}{Theorem}[section]
\newtheorem{prop}[thm]{Proposition}
\newtheorem{lem}[thm]{Lemma}
\newtheorem{clry}[thm]{Corollary}
\newtheorem{deft}[thm]{Definition}
\newtheorem{hyp}{Assumption}
\newtheorem*{helgason}{Helgason's support theorem}

\theoremstyle{definition}
\newtheorem{rem}[thm]{Remark}
\newtheorem*{acknow}{Acknowledgements}
\numberwithin{equation}{section}
\newcommand{\eps}{\varepsilon}
\renewcommand{\phi}{\varphi}
\renewcommand{\d}{\partial}
\newcommand{\re}{\mathop{\rm Re} }
\newcommand{\im}{\mathop{\rm Im}}
\newcommand{\R}{\mathbf{R}}
\newcommand{\C}{\mathbf{C}}
\newcommand{\N}{\mathbf{N}} 
\newcommand{\supp}{\mathop{\rm supp}}
\newcommand{\ch}{\mathop{\rm ch}}
\newcommand{\vararg}{\mathop{\rm vararg}}
\newcommand{\curl}{\mathop{\rm curl}}
\title[Determining a magnetic Schr\"odinger operator]{Determining a magnetic Schr\"odinger operator from partial Cauchy data}
\author[]{David Dos Santos Ferreira, Carlos E. Kenig, \\ Johannes Sj\"ostrand, Gunther Uhlmann}
\begin{abstract}
   In this paper we show, in dimension $n\geq 3$, that knowledge of the Cauchy data for the Schr{\"o}dinger equation
   in the presence of a magnetic potential, measured on possibly very small subsets of the boundary, determines uniquely the magnetic field 
   and the electric potential. We follow the general strategy of \cite{KSU} using a richer set of solutions to the Dirichlet problem
   that has been used in previous works on this problem.
\end{abstract}
\maketitle
%
%%%%%%%%%%%%%%%%%%%%%%%%%%%%%%%%%%%%%%%%%%%%%%%%%%%%%%%%%%%%%%
%%%%%%%%%%%%%%%%%%%%%%%%%%%%%%%%%%%%%%%%%%%%%%%%%%%%%%%%%%%%%%
%%                                                                              INTRODUCTION                                                                                     %%
%%%%%%%%%%%%%%%%%%%%%%%%%%%%%%%%%%%%%%%%%%%%%%%%%%%%%%%%%%%%%%
%%%%%%%%%%%%%%%%%%%%%%%%%%%%%%%%%%%%%%%%%%%%%%%%%%%%%%%%%%%%%%
%
\begin{section}{Introduction}
Let  $\Omega \subset \R^{n}$ be an open bounded set with $\mathcal{C}^{\infty}$ boundary, 
we are interested in the magnetic Schr\"odinger operator
\begin{align}
\label{intro:Schrod}
    \mathcal{L}_{A,q}(x,D)&=\sum_{j=1}^{n} (D_{j}+A_{j}(x))^{2}+q(x) \\ \nonumber
    &=D^{2}+A \cdot D + D \cdot A +A^{2}+q
\end{align}
with real magnetic potential $A=(A_{j})_{1 \leq j \leq n} \in \mathcal{C}^{2}(\bar{\Omega};\R^{n})$ and bounded electric potential $q \in L^{\infty}(\Omega)$. As usual, $D=-i\nabla$. In this paper, we always assume the dimension to be $\geq 3$. Let us introduce the
\begin{hyp}
\label{intro:eigenHyp}
     $0$ is not an eigenvalue of the magnetic Schr\"odinger operator $\mathcal{L}_{A,q}:H^{2}(\Omega) \cap H^{1}_{0}(\Omega) 
     \rightarrow L^{2}(\Omega)$.
\end{hyp}
Let $\nu$ be the unit exterior normal. Under the assumption \ref{intro:eigenHyp}, the Dirichlet problem
\begin{align}
      \left\{ \begin{array}{l} \mathcal{L}_{A,q}u=0 \\ u|_{\d\Omega}=f \in H^{\frac{1}{2}}(\d\Omega) \end{array} \right. 
 \end{align}
 has a unique solution in $H^{1}(\Omega)$, and we can introduce the Dirichlet to Neumann map (DN)
     $$ \mathcal{N}_{A,q}:H^{\frac{1}{2}}(\d \Omega) \ni f \mapsto (\d_{\nu}+iA \cdot \nu) u|_{\d \Omega} 
          \in  H^{-\frac{1}{2}}(\d\Omega) $$ 
associated to the magnetic Schr\"odinger operator $\mathcal{L}_{A,q}$ with magnetic potential defined by (\ref{intro:Schrod}). 

The inverse problem we consider in this paper is to recover information
about the magnetic and electric potential from the DN map measured on
subsets of the boundary.  As was noted in \cite{Su}, the DN map
is invariant under a gauge transformation of the magnetic potential:
it ensues from the identities
\begin{align}
\label{intro:gauge}
    e^{-i\Psi}\mathcal{L}_{A,q}e^{i\Psi}=\mathcal{L}_{A+\nabla\Psi,q}, \quad 
    e^{-i\Psi}\mathcal{N}_{A,q}e^{i\Psi}=\mathcal{N}_{A+\nabla\Psi,q},
\end{align}
that $\mathcal{N}_{A,q}=\mathcal{N}_{A+\nabla \Psi,q}$ when $\Psi \in C^{1}(\bar{\Omega}) $ is 
such that  $\Psi|_{\partial\Omega}=0$. Thus $\mathcal{N}_{A,q}$ carries information about the magnetic field%
\footnote{Here $A$ is viewed as the $1$-form $\sum_{j=1}^{n}A_{j} dx_{j}$.}
 $B=d A$. Sun showed in \cite{Su} that from this information one can determine
the magnetic field and the electric potential if the magnetic potential
is small in an appropriate class.  In \cite{NSU} the smallness assumption
was eliminated for smooth magnetic and electric potentials and for $\mathcal{C}^2$
and compactly supported magnetic potential
and $L^\infty$ electrical potential. The regularity assumption on the
magnetic potential
was improved  in  \cite{T} to $\mathcal{C}^{2/3+\epsilon}, \epsilon>0$,  and  to Dini continuous in \cite{Sa1}. Recently in \cite{Sa2}
a method was given for reconstructing the magnetic field and the electric potential
under some regularity assumptions on the magnetic potential.

All of the above mentioned results rely on constructing complex
geometrical optics solutions, with a linear phase, for the magnetic Schr\"odinger
equation. We also mention that the inverse boundary value problem is
closely related to the inverse scattering problem at a fixed energy
for the magnetic Schr\"odinger operator. The latter was
studied under various regularity assumptions on the magnetic
and electrical potentials in \cite{NK}, for small compactly supported
magnetic potential and compactly supported electric potential. This
result was extended in \cite{ER} for exponentially decaying magnetic
and electric potentials with no smallness assumption.  

In this paper we extend the main result of \cite{KSU} to the case of the magnetic Schr\"odinger
equation.  We state the precise results below.
\\

Let $x_{0} \in \R^{n} \backslash \overline{\ch(\Omega)}$ (where $\ch(\Omega)$ denotes the convex hull of $\Omega$), we define 
the front and back sides of $\d\Omega$ with respect to $x_{0}$ by
\begin{align*}
    F(x_{0})&=\{x \in \d\Omega : (x-x_{0}) \cdot \nu(x) \leq 0\} \\
    B(x_{0})&=\{x \in \d\Omega : (x-x_{0}) \cdot \nu(x) > 0\}.
\end{align*} 

\begin{thm}
\label{intro:mainThm}
     Let $\Omega$ be a simply connected open bounded set with $\mathcal{C}^{\infty}$ boundary in $\R^{n}$, $n\geq 3$,  
     let $A_{1},\,A_{2}$ be two real $\mathcal{C}^{2}$ vector fields on $\bar{\Omega}$ and $q_{1},\,q_{2}$ be two bounded  
     potentials on $\Omega$ such that the assumption \ref{intro:eigenHyp} is satisfied. 
          
     Let $x_{0} \in \R^{n} \backslash \overline{\ch(\Omega)}$, suppose that the Dirichlet to Neumann maps related to the 
     operators $\mathcal{L}_{A_{1},q_{1}}$ and $\mathcal{L}_{A_{2},q_{2}}$ coincide on part of the boundary near $x_{0}$ in the sense 
     that there exists a neighborhood $\tilde{F}$ of the front side of $\d\Omega$ with respect to $x_{0}$ such that
     \begin{align}
     \label{intro:DTNcoincide}
          \mathcal{N}_{A_{1},q_{1}}f(x)=\mathcal{N}_{A_{2},q_{2}}f(x) \quad \forall x \in \tilde{F}, 
          \quad  \forall f \in H^{\frac{1}{2}}(\d\Omega), 
     \end{align}
     then if $A_{1}$ and $A_{2}$ are viewed as $1$-forms
         $$ d A_{1}=d A_{2} \textrm{ and } q_{1}=q_{2}. $$
\end{thm}
\begin{rem}
    We only use the simple connectedness of the set $\Omega$ to deduce that the two magnetic potentials differ from a gradient,  from the equality
    $dA_{1}=dA_{2}$. If we already know that $A_{1}-A_{2}=\nabla \Psi$, we don't need the fact that $\Omega$ is simply connected,
    in particular, theorem \ref{intro:mainThm} contains theorem 1.2 of \cite{KSU}.
    
    Nevertheless theorem 1.1 in \cite{KSU} improves on this result by restricting the Dirichlet-to-Neuman maps 
    to a space of functions on the boundary with support in a small neighborhood of
    the back side $B(x_{0})$.  We have left the corresponding result in the magnetic case open.
\end{rem}
As in \cite{KSU}, we make the following definition of a strongly star shaped domain.
\begin{deft}
    An open set $\Omega$ with smooth boundary is said to be strongly star shaped with respect to $x_{1} \in \d\Omega$ if every line
    through $x_{1}$ which is not contained in the tangent hyperplane  cuts the boundary at precisely two distinct points $x_{1}$ and
     $x_{2}$ with transversal intersection at $x_{2}$. 
\end{deft}
With this definition, theorem \ref{intro:mainThm} implies the following corollary
\begin{clry}
    Under the assumptions on $\Omega$,  the magnetic potentials $A_{1}, \, A_{2}$ and the electric potentials $q_{1}, \, q_{2}$ of theorem 
    \ref{intro:mainThm},  let $x_{1} \in \d\Omega$ be a point of the boundary such that the tangent hyperplane of $\d\Omega$ at 
    $x_{1}$ only intersects $\d\Omega$ at $x_{1}$ and such that $\Omega$ is strongly star shaped with respect to $x_{1}$.
    
    Suppose that there exists a neighborhood $\tilde{F}$ of $x_{1}$ in $\d\Omega$  such that $(\ref{intro:DTNcoincide})$ holds
    then 
         $$ dA_{1}=dA_{2}  \textrm{ and }  q_{1}=q_{2}. $$
\end{clry}

We proceed as in \cite{KSU} by constructing some complex geometrical optics solutions using a Carleman estimate. The construction of these  
solutions is fairly similar to those presented in the latter paper, except for the changes due to the presence of the magnetic potential.
However, the part concerned with the recovery of the potential and the magnetic field is new. 

The plan of this article is as follows.
\tableofcontents
\end{section}
%
%%%%%%%%%%%%%%%%%%%%%%%%%%%%%%%%%%%%%%%%%%%%%%%%%%%%%%%%%%%%%%
%%%%%%%%%%%%%%%%%%%%%%%%%%%%%%%%%%%%%%%%%%%%%%%%%%%%%%%%%%%%%%
%%                                                                    CARLEMAN ESTIMATE                                                                                  %%
%%%%%%%%%%%%%%%%%%%%%%%%%%%%%%%%%%%%%%%%%%%%%%%%%%%%%%%%%%%%%%
%%%%%%%%%%%%%%%%%%%%%%%%%%%%%%%%%%%%%%%%%%%%%%%%%%%%%%%%%%%%%%
%
\begin{section}{Carleman estimate}
Our first step is to construct solutions of the magnetic Schr\"odinger equation $\mathcal{L}_{A,q}u=0$ of the form
\begin{align}
\label{Carleman:WKBsol}
     u(x,h) = e^{\frac{1}{h}(\phi+i\psi)}(a(x)+h r(x,h)) 
\end{align}
(where $\phi$ and $\psi$ are real functions) by use of the complex geometrical optics method: of course, $\psi$ and $a$ will be sought as 
solutions of respectively an eikonal equation and a transport equation. In order to be able to go from an approximate solution to an exact
solution,  one wants the conjugated operator 
    $$ e^{\frac{\phi}{h}}h^{2}\mathcal{L}_{A,q}e^{-\frac{\phi}{h}} $$
to be locally solvable in a semi-classical sense, which means its principal symbol%
\footnote{Here and throughout this article, we are using the semi-classical convention.}
\begin{align}
\label{Carleman:ConjSymb}
     p_{\phi}(x,\xi) = \xi^{2} - (\nabla \phi)^{2}+2i \nabla \phi \cdot \xi 
\end{align}
to satisfy H\"ormander's condition
    $$ \{\re p_{\phi},\im p_{\phi}\} \leq 0 \textrm{ when } p_{\phi}=0.  $$

Since we furthermore want to obtain solutions (\ref{Carleman:WKBsol}) for both the phases $\phi$ and $-\phi$,  we
will consider phases satisfying the condition
\begin{align}
\label{Carleman:Hcondition}
     \{\re p_{\phi},\im p_{\phi}\} = 0 \textrm{ when } p_{\phi}=0.  
\end{align}
\begin{deft}
    A real smooth function $\phi$ on an open set $\tilde{\Omega}$ is said to be a limiting Carleman weight if it has non-vanishing 
    gradient on $\tilde{\Omega}$ and if the symbol $(\ref{Carleman:ConjSymb})$ satisfies the condition $(\ref{Carleman:Hcondition})$  on 
    $T^{*}(\tilde{\Omega})$. This is equivalent to say that 
    \begin{align}
    \tag{\ref{Carleman:Hcondition}$'$}
          \langle \phi'' \nabla \phi,\nabla \phi \rangle + \langle\phi'' \xi,\xi\rangle=0 \textrm{ when } \xi^{2}=(\nabla \phi)^{2}
          \textrm{ and } \nabla \phi \cdot \xi=0.
    \end{align}
\end{deft}

The appropriate tool in deducing local solvability for the conjugated operator and in proving that the geometrical optics method is
 effective  (meaning that indeed one gains one power of $h$ in the former asymptotics) is a Carleman estimate. The goal
of this section is to prove such an estimate.

In this section, $\Omega$ is as in the introduction and $\tilde{\Omega}$ will denote an open set $\tilde{\Omega} \Supset \Omega$.
We will use the following notations
     $$ (u|v)=\int_{\Omega} u(x) \bar{v}(x) \, dx, \quad  (u|v)_{\d\Omega}=\int_{\d\Omega} u(x) \bar{v}(x) 
          \, d\sigma(x), $$
and $\|u\|=\sqrt{(u|u)}$ denotes the $L^{2}$ norm on $\Omega$. We say that the estimate 
    $$ F(u,h) \lesssim G(u,h) $$
holds for all $u \in X$ (where $X$ is a function space, such as $L^{2}(\Omega)$)  and for $h$ small if there exist constants $C>0$
and $h_{0}>0$ (possibly depending on $q$ and $A$) such that for all $0\leq h \leq h_{0}$ and for all $u \in X$, the inequality 
$F(u,h) \leq C G(u,h)$ is satisfied.

We will make extensive use of the Green formula for the magnetic Schr\"odinger operator $\mathcal{L}_{A,q}$, which for sake of convenience, we
state as a lemma.
\begin{lem}
     Let $A$ be a real $\mathcal{C}^{1}$ vector field on $\bar{\Omega}$ and $q \in L^{\infty}(\Omega)$ 
     then we have the magnetic Green formula
     \begin{multline}
     \label{Carleman:Green}
          (\mathcal{L}_{A,q}u|v)_{\Omega}-(u|\mathcal{L}_{A,\bar{q}}v)_{\Omega} =\\  \big(u|(\d_{\nu}+i\nu \cdot A)v\big)_{\d\Omega}
          -\big((\d_{\nu}+i\nu \cdot A)u|v\big)_{\d \Omega}
     \end{multline}  
     for all $u, \, v \in H^{1}(\Omega)$ such that $\Delta u, \, \Delta v \in L^{2}(\Omega)$. 
\end{lem}
\begin{proof}
     Integrating by parts, we have
      \begin{align}
      \label{Carleman:IbyP}
           (\mathcal{L}_{A,q}u|v)_{\Omega}&=(\nabla u|\nabla v)_{\Omega}+(Du|Av)_{\Omega}+(Au|Dv)_{\Omega}
           \\ \nonumber &\quad +\big((q+A^{2})u|v\big)_{\Omega}-\big((\d_{\nu}+i\nu \cdot A)u|v\big)_{\d \Omega}
      \end{align}
      and by permuting $u$ and $v$, replacing $q$ by $\bar{q}$, and taking the complex conjugate of the former, we get 
       \begin{align*}
           (\mathcal{L}_{A,\bar{q}}u|v)_{\Omega}&=(\nabla u|\nabla v)_{\Omega}+(Du,Av)_{\Omega}+(Au|Dv)_{\Omega}
           \\ \nonumber &\quad +\big((q+A^{2})u|v\big)_{\Omega}-\big(u|(\d_{\nu}+i\nu \cdot A)v\big)_{\d \Omega}.
       \end{align*}
       Subtracting the former to (\ref{Carleman:IbyP}), we end up with (\ref{Carleman:Green}).
\end{proof}
If $\phi$ is a limiting Carleman weight, we define
   $$  \d\Omega_{\pm}=\{x \in \d\Omega : \pm \d_{\nu}\phi \geq 0\}. $$
\begin{prop}
\label{Carleman:Carlprop}
   Let $\phi$ be a $\mathcal{C}^{\infty}$ limiting Carleman weight on $\tilde{\Omega}$, 
   let $A$ be a $\mathcal{C}^{1}$ vector field on $\bar{\Omega}$ and $q \in L^{\infty}(\Omega)$, the Carleman estimate
   \begin{multline}
   \label{Carleman:CarlEstBdary}
        -h(\d_{\nu}\phi \, e^{\frac{\phi}{h}}\d_{\nu}u|e^{\frac{\phi}{h}}\d_{\nu}u)_{\d\Omega_{-}}+
        \|e^{\frac{\phi}{h}}u\|^{2} +\|e^{\frac{\phi}{h}}h\nabla u\|^{2}\\
        \lesssim h^{2} \|e^{\frac{\phi}{h}}\mathcal{L}_{A,q}u\|^{2}+
        h(\d_{\nu}\phi \, e^{\frac{\phi}{h}}\d_{\nu}u|e^{\frac{\phi}{h}}\d_{\nu}u)_{\d\Omega_{+}}
   \end{multline}
   holds for all $u \in \mathcal{C}^{\infty}(\bar{\Omega}) \cap H^{1}_{0}(\Omega)$ and $h$ small.
\end{prop}
In particular, when $u \in \mathcal{C}^{\infty}_{0}(\Omega)$, we have the Carleman estimate
\begin{align}
\tag{\ref{Carleman:CarlEstBdary}$'$}
        \|e^{\frac{\phi}{h}}u\| +h\|e^{\frac{\phi}{h}}\nabla u\| \lesssim h \|e^{\frac{\phi}{h}}\mathcal{L}_{A,q}u\|.
\end{align}
\begin{proof}
   Taking $v=e^{\frac{\phi}{h}}u$, it is equivalent to prove the following \textit{ a priori} estimate
    \begin{multline}
    \label{Carleman:Apriori}
         -h (\d_{\nu}\phi \, \d_{\nu}v|\d_{\nu}v)_{\d\Omega_{-}} + \|v\|^{2}+\|h\nabla v\| ^{2} \\
         \lesssim \frac{1}{h^{2}} \|(e^{\frac{\phi}{h}}h^{2}\mathcal{L}_{A,q}e^{-\frac{\phi}{h}})v\|^{2}
        +h (\d_{\nu}\phi \, \d_{\nu}v|\d_{\nu}v)_{\d\Omega_{+}}
    \end{multline}
    since $\|he^{\frac{\phi}{h}}\nabla u\|  \lesssim \|v\|+\|h\nabla v\| $ and $(e^{\frac{\phi}{h}}\d_{\nu}u)|_{\d\Omega}=
    \d_{\nu}v|_{\d \Omega}$. Conjugating the magnetic Schr\"odinger operator by the exponential weight gives rise to the following operator
    \begin{align}
         e^{\frac{\phi}{h}}h^{2}\mathcal{L}_{A,q}e^{-\frac{\phi}{h}} = P+iQ+R+h^{2}(q+A^{2})
    \end{align}
    where $P$ and $Q$ are the self-adjoint operators
    \begin{align*}
        P&= h^{2}D^{2}-(\nabla \phi)^{2},  \\
        Q&=\nabla \phi \cdot hD + hD \cdot \nabla \phi, \\
        \textrm{and} \quad R&=h (A \cdot hD +hD \cdot A)+2i h A \cdot \nabla \phi.
    \end{align*}
    Our first remark concerns the fact that we may neglect the term $h^{2}(q+A^{2})$ since the right-hand side of
    (\ref{Carleman:Apriori}) may be perturbed by a term bounded by $h^{2} \|v\|^{2}$, which may be absorbed by  
    the left-hand side if $h$ is small enough. Omitting the term $q+A^{2}$ gives rise to such an error.
    Hence we will prove the \textit{a priori} estimate for the operator $P+iQ+R$.  The same is not true
    of the term $R$ because errors of order $\|v\|^{2}+\|h\nabla v\|^{2}$ may not be absorbed into the left hand-side.
    
    Note that if $p$ and $q$ denote the principal symbol respectively of $P$ and $Q$, the fact that $\phi$ is a limiting Carleman 
    weight means that 
         $$ \{p,q\}=0 \textrm{ when } p+iq=0. $$ 
    This condition is not enough to obtain an \textit{a priori} estimate for $P+iQ$, one needs to have a positive Poisson bracket. Our
    first step is to remedy to this by using a classical convexity argument. Consider the modified Carleman weight 
         $$ \tilde{\phi}= \phi + h \frac{\phi^{2}}{2 \eps} $$ 
    where $\eps$ is a suitable small parameter to be chosen independent of $h$, and denote by $\tilde{p}$ and $\tilde{q}$ 
    the corresponding symbols, and by $\tilde{P},\tilde{Q},\tilde{R}$ the corresponding operators,
    when $\phi$ has been replaced by $\tilde{\phi}$. Then, we have
        $$ \nabla \tilde{\phi}=(1+\frac{h}{\eps} \phi)\nabla \phi, \quad \tilde{\phi}''=(1+\frac{h}{\eps} \phi) \phi''
             +\frac{h}{\eps} \nabla \phi \otimes \nabla \phi $$
    therefore when $\xi^{2}=(\nabla \tilde{\phi})^{2}$ and $\nabla \tilde{\phi} \cdot \xi=0$, we have 
    \begin{align}
    \label{Carleman:convexBracket}
         \{\tilde{p},\tilde{q}\}&=4\langle \tilde{\phi}'' \xi,\xi \rangle +4\langle \tilde{\phi}'' \nabla \tilde{\phi},\nabla \tilde{\phi} \rangle 
         \\ \nonumber &= \frac{4h}{\eps} (1+\frac{h}{\eps}\phi)^{2}(\nabla \phi)^{4}>0
    \end{align}
    since $\phi$ is a limiting Carleman weight.  Furthermore, if we restrict ourselves to the hyperplane $V_{x}$ orthogonal to $\nabla \phi$, 
    we get 
        $$  \{\tilde{p},\tilde{q}\}(x,\cdot)|_{V_{x}}=\frac{4h}{\eps} (1+\frac{h}{\eps}\phi)^{2}(\nabla \phi)^{4} 
              + a(x)(\xi^{2}-(\nabla \tilde{\phi})^{2}) $$
    with $a(x)=4h(\nabla \tilde{\phi})^{2}/\eps-4\langle \tilde{\phi}'' \nabla \tilde{\phi},\nabla \tilde{\phi} \rangle 
    / (\nabla \tilde{\phi})^{2}$, and since this bracket
    is a quadratic polynomial with no linear part, this implies that there exists a linear form $b(x,\xi)$ in $\xi$ such that
        $$  \{\tilde{p},\tilde{q}\}=\frac{4h}{\eps} (1+\frac{h}{\eps}\phi)^{2}(\nabla \phi)^{4}+a(x) \tilde{p}+b(x,\xi) \tilde{q}. $$   
    
    This computation implies on the operator level that
    \begin{align*}
        i[\tilde{P},\tilde{Q}]&=\frac{4h^{2}}{\eps} (1+\frac{h}{\eps}\phi)^{2}(\nabla \phi)^{4}+\frac{h}{2}(a \tilde{P}+\tilde{P} a)
        \\  &\quad +\frac{h}{2}(b^{w} \tilde{Q}+\tilde{Q}b^{w})+h^{3} c(x)
    \end{align*}
    where  the first order differential operator $b^{ w}$ is the semi-classical Weyl quantization%
    \footnote{The absence of $h^{2}$ term is due to the use of the Weyl quantization.} of $b$.  In fact, the positivity of the bracket
    (\ref{Carleman:convexBracket}) essentially induces the positivity of the commutator $i[\tilde{P},\tilde{Q}]$
    \begin{align}
    \label{Carleman:commut}
         i([\tilde{P},\tilde{Q}]v|v) &= \underbrace{\frac{4h^{2}}{\eps} (1+\frac{h}{\eps}\phi)^{2} \|(\nabla \phi)^{2}v\|^{2}}_{>0}
         \\ \nonumber &\quad    + h \re (a\tilde{P}v|v) + h \re(\tilde{Q}v| b^{w}v)+h^{3}(cv|v) 
    \end{align}
    (recall that $v|_{\d \Omega}=0$ which explains why  there are no boundary terms). The former fact will be enough to obtain the \textit{a priori} 
    estimate on $\tilde{P}+i\tilde{Q}$. 
    
    Our last observation is that
    \begin{align*}
         \|h \nabla v\|^{2} = (\tilde{P}v|v) + \| \sqrt{(\nabla \phi)^{2}} \, v\|^{2}
    \end{align*}
    leading to 
    \begin{align}
    \label{Carleman:estDer}
         \|h \nabla v\|^{2} \lesssim  \|\tilde{P}v\|^{2}+\|v\|^{2}.
    \end{align}
        
    Now, we turn to the proof of the estimate. We have
    \begin{align*}
         \|(\tilde{P}+i\tilde{Q})v\|^{2}=\|\tilde{P}v\|^{2} + \|\tilde{Q}v\|^{2}+i(\tilde{Q}v|\tilde{P}v)-i(\tilde{P}v|\tilde{Q}v)
    \end{align*}
    and the magnetic Green formula (\ref{Carleman:Green})  (used in the straightforward case with no potential  
    $\tilde{P}=h^{2}\mathcal{L}_{0,-(\nabla \tilde{\phi})^{2}/h^{2}}$), together with the fact that $v|_{\d \Omega}=0$, gives
    \begin{align*}
         (\tilde{Q}v|\tilde{P}v)&=(\tilde{P}\tilde{Q}v|v)-h^{2}(\tilde{Q} v| \d_{\nu}v)_{\d\Omega} \\ &= 
         (\tilde{P}\tilde{Q}v|v)+2ih^{3}(\d_{\nu}\tilde{\phi} \,  \d_{\nu}v| \d_{\nu}v)_{\d\Omega} 
    \end{align*}
    and similarly, since $\tilde{Q}$ is first order, we get 
         $$ (\tilde{P}v,\tilde{Q}v)=(\tilde{Q}\tilde{P}v,v). $$
    Therefore we have 
     \begin{align*}
         \|(\tilde{P}+i\tilde{Q})v\|^{2}&=\|\tilde{P}v\|^{2} + \|\tilde{Q}v\|^{2}+i([\tilde{P},\tilde{Q}]v|v)
         \\ &\quad -2h^{3}( \d_{\nu}v| \d_{\nu}\tilde{\phi} \, \d_{\nu}v)_{\d\Omega}
     \end{align*}
     and using (\ref{Carleman:commut}), we get 
     \begin{multline*}
         \|(\tilde{P}+i\tilde{Q})v\|^{2}+2h^{3}( \d_{\nu}v| \d_{\nu}\tilde{\phi} \, \d_{\nu}v)_{\d\Omega}
         \geq \|\tilde{P}v\|^{2} + \|\tilde{Q}v\|^{2}+ \\ +\frac{C_{1}h^{2}}{\eps} \|v\|^{2} 
         -\underbrace{(Ch\|\tilde{P}v\| \, \|v\|+Ch\|\tilde{Q}v\| \, \|h\nabla v\|)}_{\leq \frac{1}{2}\|\tilde{P}v\|^{2}
          +\frac{1}{2}\|\tilde{Q}v\|^{2}+\frac{C_{2}^{2}h^{2}}{2}\big(\|v\|^{2} +\|h\nabla v\|^{2} \big)}
     \end{multline*}
     which combined with (\ref{Carleman:estDer}), gives when $\eps$ is small enough   
     \begin{multline*}
         \|(\tilde{P}+i\tilde{Q})v\|^{2}+2h^{3}( \d_{\nu}v| \d_{\nu}\tilde{\phi} \, \d_{\nu}v)_{\d\Omega}
         \\ \gtrsim (1-\mathcal{O}(\eps^{-1}h^{2}) \|\tilde{P}v\|^{2} + \|\tilde{Q}v\|^{2}+\frac{h^{2}}{\eps} 
         \big( \|v\|^{2}+\|h\nabla v\|^{2} \big).
     \end{multline*}
     Thus, taking $h$ and $\eps$ small enough
     \begin{align}
      \label{Carleman:beforeCarl}
         \|(\tilde{P}+i\tilde{Q})v\|^{2}+2h^{3}( \d_{\nu}v| \d_{\nu}\tilde{\phi} \, \d_{\nu}v)_{\d\Omega}
         \gtrsim \frac{h^{2}}{\eps}\big( \|v\|^{2}+\|h\nabla v\|^{2} \big).
     \end{align}
     
     The last part%
     \footnote{\label{modifCarl}This is the main difference with respect to the proof of the Carleman estimate in \cite{KSU}.} 
     of the proof is concerned with the additional term $\tilde{R}v$ due to the magnetic potential; from the former
     inequality we deduce
     \begin{multline*}
         \|(\tilde{P}+i\tilde{Q}+\tilde{R})v\|^{2}+h^{3}( \d_{\nu}v| \d_{\nu}\tilde{\phi} \, \d_{\nu}v)_{\d\Omega}
         \\ \gtrsim \frac{h^{2}}{\eps}\big( \|v\|^{2}+\|h\nabla v\|^{2}-\mathcal{O}(\eps) \|h^{-1}\tilde{R} v\|^{2}\big)
     \end{multline*}
     and using the fact that $\|h^{-1}\tilde{R} v\|^{2} \lesssim \|v\|^{2}+\|h\nabla v\|^{2}$, we obtain
      \begin{align*}
         \|(\tilde{P}+i\tilde{Q}+\tilde{R})v\|^{2}+h^{3}( \d_{\nu}v| \d_{\nu}\tilde{\phi} \, \d_{\nu}v)_{\d\Omega}
         \gtrsim \frac{h^{2}}{\eps}\big( \|v\|^{2}+\|h\nabla v\|^{2}\big)
     \end{align*}
     if $\eps$ is chosen small enough. Finally, with $v=e^{\frac{\phi^{2}}{2\eps}}e^{\frac{\phi}{h}}u$, we get
     \begin{multline*}
         \|e^{\frac{\phi^{2}}{2\eps}}e^{\frac{\phi}{h}}u\|^{2} + \|h \nabla e^{\frac{\phi^{2}}{2\eps}}e^{\frac{\phi}{h}}u\|^{2} 
         \lesssim \frac{1}{h^{2}} \|e^{\frac{\phi^{2}}{2\eps}}e^{\frac{\phi}{h}}h^{2}\mathcal{L}_{A,q}u\|^{2}
         \\ +h\big(e^{\frac{\phi^{2}}{\eps}} \d_{\nu}(e^{\frac{\phi}{h}}u)| 
         \d_{\nu}\tilde{\phi} \, \d_{\nu}(e^{\frac{\phi}{h}}u)\big)_{\d\Omega}
     \end{multline*}
     this gives the Carleman estimate (\ref{Carleman:CarlEstBdary}) since 
          $$ 1 \leq e^{\phi^{2}/2\eps} \leq C, \quad 
               \frac{1}{2}  \leq \frac{\d_{\nu}\tilde{\phi}}{\d_{\nu}\phi}=1+\frac{h}{\eps}\phi \leq \frac{3}{2} $$ 
     on $\bar{\Omega}$ for all $h$ small enough.
\end{proof}

We denote by $H^{1}_{\rm scl}(\Omega)$ the semi-classical Sobolev space of order $1$ on $\Omega$ with associated norm
    $$ \|u\|^2_{H^{1}_{\rm scl}(\Omega)}=\| h\nabla u\|^{2}+\|u\|^{2} $$
and by $H^{s}_{\rm scl}(\R^{n})$ the semi-classical Sobolev space on $\R^{n}$ with associated norm
    $$ \|u\|^2_{H^{s}_{\rm scl}(\R^{n})}=\|\langle hD \rangle^{s}u\|^{2}_{L^{2}(\R^{n})}
          =\int (1+h^{2}\xi^{2})^{s}|\hat{u}(\xi)|^{2} \, d\xi. $$

Changing $\phi$ into $-\phi$, we may rewrite the Carleman estimate in the following convenient way 
\begin{multline}
\label{Carleman:CarlEst}
        \sqrt{h} \|\sqrt{\d_{\nu}\phi} \,  e^{-\frac{\phi}{h}}\d_{\nu}u\|_{L^{2}(\d\Omega_{+})}
        +\|e^{-\frac{\phi}{h}} u\|_{H^{1}_{\rm scl}(\Omega)} \\ \lesssim h \|e^{-\frac{\phi}{h}}\mathcal{L}_{A,q}u\|
        +\sqrt{h} \|\sqrt{-\d_{\nu}\phi} \,  e^{-\frac{\phi}{h}}\d_{\nu}u\|_{L^{2}(\d\Omega_{-})}.
\end{multline}
By regularization, this estimate is still valid for $u \in H^{2}(\Omega) \cap H^{1}_{0}(\Omega)$.
A similar Carleman estimate gives the following solvability result:
\begin{prop}
\label{Carleman:solvability}
    Let $\phi$ be a limiting Carleman weight on $\tilde{\Omega}$, let $A$ be a $\mathcal{C}^{1}$ vector field on $\bar{\Omega}$
    and $q \in L^{\infty}(\Omega)$. There exists $h_{0}$ such that for all $0 \leq h \leq h_{0}$ and for all $w \in L^{2}(\Omega)$,
    there exists 
    $u \in H^{1}(\Omega) $ such that 
        $$ h^{2}\mathcal{L}_{A,q} (e^{\frac{\phi}{h}}u)= e^{\frac{\phi}{h}}w \quad \textrm{and} \quad
             h \|u\|_{H^{1}_{\rm scl}(\Omega)}  \lesssim \|w\| . $$
\end{prop}
\begin{proof}
     We need the following  Carleman estimate 
     \begin{align}
     \label{Carleman:shiftedCarl}
          \|v\| \lesssim h \|e^{-\frac{\phi}{h}}\mathcal{L}_{A,q}e^{\frac{\phi}{h}}v\| _{H^{-1}_{\rm scl}(\R^{n})}, \quad 
               \forall v \in \mathcal{C}^{\infty}_{0}(\Omega). 
     \end{align}  
     Let $\Omega \Subset \hat{\Omega} \subset \tilde{\Omega}$, assume that we have extended $A$ to a $\mathcal{C}^{1}$ vector
     field on $\hat{\Omega}$  and $q$ to a $L^{\infty}$ function on $\hat{\Omega}$. Let $\chi \in \mathcal{C}^{\infty}_{0}(\hat{\Omega})$
     equal 1 on $\Omega$.
     With the notations used in the proof of proposition \ref{Carleman:Carlprop} we have
         $$ \langle h D \rangle^{-1} (\tilde{P}+i\tilde{Q}) \langle h D \rangle = \tilde{P}+i\tilde{Q}+hR_{1} $$
     where $R_{1}$ is a semi-classical pseudo-differential operator of order 1,
     therefore from estimate (\ref{Carleman:beforeCarl}) we deduce
     \begin{align*}
         \|(\tilde{P}+i\tilde{Q})\langle h D \rangle v\|^{2}_{H^{-1}_{\rm scl}(\R^{n})} &\gtrsim \frac{h^{2}}{\eps}
         (\| v \|^{2}_{H^{1}_{\rm scl}(\R^{n})}-\mathcal{O}(\eps)\|R_{1}v\|^{2}) \\ &\gtrsim \frac{h^{2}}{\eps}
         \| v \|^{2}_{H^{1}_{\rm scl}(\R^{n})}
     \end{align*}
     for any $v \in \mathcal{C}^{\infty}_{0}(\hat{\Omega})$, if $h$ and $\eps$ are small enough. Besides, we have 
     \begin{align*}
         \|(\tilde{R}+h^{2}(q+A^{2})) v\|_{H^{-1}_{\rm scl}(\R^{n})} \lesssim h \|v\|
     \end{align*}
     therefore if $\eps$ is small enough, we have
       $$ \|\langle h D \rangle v\| \lesssim h 
             \|e^{-\frac{\tilde{\phi}}{h}}\mathcal{L}_{A,q}e^{\frac{\tilde{\phi}}{h}}\langle h D \rangle v\| _{H^{-1}_{\rm scl}(\R^{n})} $$ 
      for any $v \in \mathcal{C}^{\infty}_{0}(\hat{\Omega})$. Hence if $u \in \mathcal{C}^{\infty}_{0}(\Omega)$,
      taking $v=\chi \langle h D \rangle^{-1}u \in \mathcal{C}^{\infty}(\hat{\Omega})$ in the former estimate,  and using the fact that
           $$ \| (1-\chi) \langle hD \rangle^{-1}u\|_{H^{s}_{\rm scl}}=\mathcal{O}(h^{\infty}) \|u\| $$      
     we obtain
          $$ \|u\| \lesssim h \|e^{-\frac{\tilde{\phi}}{h}}\mathcal{L}_{A,q}e^{\frac{\tilde{\phi}}{h}}u\| _{H^{-1}_{\rm scl}(\R^{n})}, \quad 
               \forall u \in \mathcal{C}^{\infty}_{0}(\Omega). $$     
     This gives the Carleman estimate (\ref{Carleman:shiftedCarl}) since $e^{\tilde{\phi}/h}=e^{\phi^{2}/\eps}e^{\phi/h}$.
     Classical arguments involving the Hahn-Banach theorem give the solvability result.  
\end{proof}
\end{section}
%
%%%%%%%%%%%%%%%%%%%%%%%%%%%%%%%%%%%%%%%%%%%%%%%%%%%%%%%%%%%%%%
%%%%%%%%%%%%%%%%%%%%%%%%%%%%%%%%%%%%%%%%%%%%%%%%%%%%%%%%%%%%%%
%%                                                                             WKB SOLUTIONS                                                                                   %%
%%%%%%%%%%%%%%%%%%%%%%%%%%%%%%%%%%%%%%%%%%%%%%%%%%%%%%%%%%%%%%
%%%%%%%%%%%%%%%%%%%%%%%%%%%%%%%%%%%%%%%%%%%%%%%%%%%%%%%%%%%%%%
%
\begin{section}{Construction of solutions by complex geometrical optics}
\label{WKB}
The goal of this section is to construct solutions of the magnetic Schr\"odinger equation of the form (\ref{Carleman:WKBsol}).  To do so we 
take $\psi$ to be a solution of the eikonal equation
      $$ p(x,\nabla \psi(x))+iq(x,\nabla \psi(x)) = 0 $$
such solutions exist since $\{p,q\}=0$ when $p=q=0$.  More precisely, the eikonal equation reads
\begin{align}
\label{WKB:eikonal}
    (\nabla \psi)^{2}=(\nabla \phi)^{2}, \quad \nabla \phi \cdot \nabla \psi =0.
\end{align}
In fact, in the remainder of this article, we fix the limiting Carleman weight to be 
\begin{align}
\label{WKB:phi}
    \phi(x)=\frac{1}{2} \log(x-x_{0})^{2}.  
\end{align}

For such a choice of $\phi$,  the second part of the eikonal equation is merely the fact that $\psi$
is a function of the angular variable $(x-x_{0})/|x-x_{0}|$ and we can actually give an explicit solution of the eikonal equation
\begin{align}
\label{WKB:psi}
    \psi(x)&=\frac{\pi}{2}-\arctan \frac{\omega \cdot (x-x_{0})}{\sqrt{(x-x_{0})^{2}-(\omega \cdot (x-x_{0}))^{2}}}
     \\\nonumber &=d_{S^{n-1}}\Big(\frac{x-x_{0}}{|x-x_{0}|},\omega\Big)
\end{align}
where $\omega \in S^{n-1}$.  Let us be more precise about the set where $\omega$ may vary, keeping in mind  that
we want this function to be smooth  --- in particular, we have to ensure that $\omega \neq (x-x_{0})/|x-x_{0}|$ whenever
$x \in \bar{\Omega}$. 

For that purpose,  let $r_{0}>0$ be large enough so that $\bar{\Omega} \subset B(x_{0},r_{0})$, 
let $H$ denote a hyperplane separating $x_{0}$ and $\ch(\Omega)$, and $H^{+}$ the open half space containing $\bar{\Omega}$ (and therefore $x_{0} \notin H^{+}$), we set
     $$ \Gamma = \{\theta \in S^{n-1}: x_{0}+r_{0} \theta \in H^{+}\}  $$
and $\check{\Gamma}$ the image of $\Gamma$ under the antipodal application. 
Let $\omega_{0} \in S^{n-1}\backslash (\Gamma \cup \check{\Gamma})$ and $\Gamma_{0}$ be a neighborhood of $\omega_{0}$
in $S^{n-1}\backslash (\Gamma \cup \check{\Gamma})$, then the distance $\Gamma \times \Gamma_{0} \ni (\theta,\omega) \rightarrow 
d_{S^{n-1}}(\theta,\omega)$ is a $\mathcal{C}^{\infty}$ function. Moreover, $\bar{\Omega} \subset \tilde{\Omega}=
x_{0}+\R_{+} \Gamma$ hence we have  $(x-x_{0})/|x-x_{0}| \in \Gamma$ for all $x$ in the open neighborhood $\tilde{\Omega}$
of $\bar{\Omega}$,  thus $\psi$  depends smoothly on the variables $(x,\omega)$ on $ \tilde{\Omega} \times \Gamma_{0}$.   
\begin{rem}
\label{WKB:remzvar}
     Suppose that $x_0=0$ and $\omega=(1,0,\dots,0)$, which we can always assume by doing a translation and a rotation. Notice
      that  by considering the complex variable $z=x_{1}+i|x'| \in \C$ (with $x=(x_{1},x') \in \R \times \R^{n-1}$), we have
          $$ \phi = \log |z|=\re \log z, \quad \psi = \frac{\pi}{2}-\arctan \frac{\re  z}{\im z}=\im \log z $$
     when $\im z>0$ (note that $\psi=\arctan(\im z/\re z)$ on the first quadrant $\re z>0, \, \im z>0$) hence $\phi+i\psi = \log z$.
\end{rem}
With such $\phi$ and $\psi$, we have
\begin{align*}
     h^{2}\mathcal{L}_{A,q}e^{\frac{1}{h}(\phi+i\psi)}&=e^{\frac{1}{h}(\phi+i\psi)}\big(h(D+A)\cdot(\nabla \psi-i\nabla \phi)
     \\ \nonumber &\quad+h(\nabla \psi-i\nabla \phi)\cdot (D+A)+h^{2}\mathcal{L}_{A,q}\big)
\end{align*}
thus we will have
\begin{align*}
  h^{2} \mathcal{L}_{A,q}\big(e^{\frac{1}{h}(\phi+i\psi)}a\big)=\mathcal{O}(h^{2})e^{\frac{1}{h}(\phi+i\psi)}
\end{align*}
if $a$ is a $\mathcal{C}^{2}$ solution of the first transport equation,  given by
\begin{align*}
   \big((D+A)\cdot(\nabla \psi-i\nabla \phi)+(\nabla \psi-i\nabla \phi)\cdot (D+A)\big)a=0.
\end{align*}
We write the latter as a vector field equation
\begin{align}
\label{WKB:TransportEq}
   (\nabla \psi-i\nabla \phi) \cdot Da+(\nabla \psi-i\nabla \phi)\cdot A a +\frac{1}{2i} (\Delta \psi-i\Delta \phi)a  =0.
\end{align}
We seek $a$ under exponential form $a=e^{\Phi}$, which means finding $\Phi$ solution of
\begin{align}
\label{WKB:TransportPhi}
    (\nabla \phi+i\nabla \psi) \cdot \nabla \Phi+i(\nabla \phi+i\nabla \psi)\cdot A +\frac{1}{2} \Delta (\phi+i \psi) =0.
\end{align}
on $\Omega$. The function $\Phi$ has $\mathcal{C}^{2}$ regularity since the magnetic potential $A$ is $\mathcal{C}^{2}$.
\begin{rem}
\label{WKB:remzvarbis}
     Considering the complex variable $z=x_{1}+i|x'|$ as in remark \ref{WKB:remzvar} with $\phi+i\psi=\log z$, 
     we may seek $\Phi$ as a solution of the following Cauchy-Riemann equation in the $z$ variable
         $$ \frac{\d \Phi}{\d \bar{z}}-\frac{(n-2)}{2(z-\bar{z})}+\frac{1}{2} A \cdot (e_{1}+ie_{r}) =0 $$
     where $e_{r}=(0,\theta)$ is the unit vector pointing in the direction of the $r$-axis.
     Indeed, if we denote by $(x_{1},r,\theta) \in \R \times \R_{+} \times S^{n-2}$ a choice of cylindrical coordinates on $\R^{n}$
     and $z=x_{1}+ir$, we have
     \begin{align*}
         \nabla(\phi + i\psi) \cdot \nabla &= \frac{\d \log z}{\d x_{1}} \, \frac{\d}{\d x_{1}} + \frac{\d \log z}{\d r} \, \frac{\d}{\d r}
         =\frac{2}{z}\frac{\d}{\d \bar z}
     \\  \textrm{ and } \quad 
          \Delta (\phi+i\psi) &= \Big(\frac{\d^{2}}{\d x_{1}^{2}}+\frac{\d^{2}}{\d r^{2}}+\frac{(n-2)}{r}\frac{\d}{\d r}+
          \frac{1}{r^{2}} \Delta_{S^{n-2}}\Big)(\phi+i\psi) \\
          \nonumber &= \frac{n-2}{r} \frac{\d}{\d r} \log z = \frac{(n-2)i}{r z}. 
     \end{align*}
\end{rem}
\begin{rem}
\label{WKB:remhol}
     Note that the set of solutions of  (\ref{WKB:TransportEq}) is invariant under the multiplication by a function $g$ satisfying
           $$ (\nabla \phi+i\nabla \psi)\cdot \nabla g =0. $$    
     In the setting of remark \ref{WKB:remzvar}, this condition reads
          $$ \frac{\d g}{\d \bar{z}} = 0 $$
     on $\Omega$, i.e. $g$ is a holomorphic function of $z=x_{1}+i|x'|$.
\end{rem}

Having chosen the phase $\phi+i\psi$ and the amplitude $e^{\Phi}$, we obtain an approximate solution of the magnetic Schr\"odinger equation
    $$ h^{2}\mathcal{L}_{A,q}(e^{\frac{1}{h}(\phi+i\psi)}e^{\Phi}) = e^{\frac{1}{h}(\phi+i\psi)}h^{2}
         \mathcal{L}_{A,q}e^{\Phi}=\mathcal{O}(h^{2}) e^{\frac{\phi}{h}} $$
(recall that $\Phi$ is $\mathcal{C}^{2}$) which we can transform into an exact solution thanks to proposition \ref{Carleman:solvability}; 
there exists $r(x,h) \in H^{1}(\Omega)$ such that 
    $$ h^{2}e^{\frac{1}{h}(\phi+i\psi)}\mathcal{L}_{A,q}r(x,h)=-e^{\frac{1}{h}(\phi+i\psi)}h \mathcal{L}_{A,q}e^{\Phi}$$
and $\|r\|_{H^{1}_{\rm scl}(\Omega)} \lesssim \|\mathcal{L}_{A,q}e^{\Phi}\|$.

We sum up the result of this section in the following lemma.
\begin{lem}
\label{WKB:WKB}
    Let $x_{0} \in \R^{n} \backslash \overline{\ch \Omega}$, there exists $h_{0}>0$ and $r$ such that 
    $\|r\|_{H^{1}_{\rm scl}(\Omega)} =\mathcal{O}(1)$ and 
        $$ u(x,h)=e^{\frac{1}{h}(\phi+i\psi)}(e^{\Phi(x)}+hr(x,h)) $$
    is a solution of the equation $\mathcal{L}_{A,q}u=0$, when $h\leq h_{0}$, 
    and $\phi$ is the limiting Carleman weight $(\ref{WKB:phi})$, $\psi$ is
    given by $(\ref{WKB:psi})$ and $\Phi$ is a solution of the Cauchy-Riemann equation $(\ref{WKB:TransportPhi})$.
\end{lem}
Note that with $\phi$ as in (\ref{WKB:phi}) the parts of the boundary $\d\Omega_{\pm}$ delimited by the sign of $\d_{\nu}\phi$ 
correspond to the front and back sides of the boundary
    $$ \d\Omega_{-}=F(x_{0}), \quad \d \Omega_{+}=\overline{B(x_{0})}. $$  
\end{section}
%
%%%%%%%%%%%%%%%%%%%%%%%%%%%%%%%%%%%%%%%%%%%%%%%%%%%%%%%%%%%%%%
%%%%%%%%%%%%%%%%%%%%%%%%%%%%%%%%%%%%%%%%%%%%%%%%%%%%%%%%%%%%%%
%%                                          RECOVERING THE MAGNETIC FIELD   		      		                                               %%
%%%%%%%%%%%%%%%%%%%%%%%%%%%%%%%%%%%%%%%%%%%%%%%%%%%%%%%%%%%%%%
%%%%%%%%%%%%%%%%%%%%%%%%%%%%%%%%%%%%%%%%%%%%%%%%%%%%%%%%%%%%%%
%
\begin{section}{Towards recovering the magnetic field}
\label{section:recovery}
Let $x_{0} \in \R^{n} \backslash \Omega$, suppose that the assumptions of theorem \ref{intro:mainThm} are fulfilled and consider
    $$ F_{\eps}=\{x \in \d\Omega : (x-x_{0}) \cdot \nu(x) <\eps |x-x_{0}| ^{2}\} \supset F(x_{0}) $$
with $\eps>0$ small enough so that  $F_{\eps} \subset \tilde{F}$, therefore satisfying
\begin{equation}
\label{recovery:hyp}
     \mathcal{N}_{A_{1},q_{1}}f(x)=\mathcal{N}_{A_{2},q_{2}}f(x), \quad \forall x \in F_{\eps},
     \quad \forall f \in H^{\frac{1}{2}}(\d\Omega). 
\end{equation}

We may assume without loss of generality that the  normal components of $A_{1}$ and $A_{2}$ are equal on the boundary
\begin{align}
\label{recovery:normComp}
    A_{1} \cdot \nu = A_{2} \cdot \nu \textrm{ on } \d\Omega
\end{align}
since we can do a gauge transformation in the magnetic potential
    $$  \mathcal{N}_{A,q}=\mathcal{N}_{A+\nabla \Psi,q} $$
(see (\ref{intro:gauge}) in the introduction) with $\Psi \in C^{3}(\bar{\Omega}) $ such that $\Psi |_{\d \Omega}=0$
and $\d_{\nu} \Psi$ is a prescribed $\mathcal{C}^{2}$ function on the boundary%
\footnote{After the use of a partition of unity and a transfer to $\{x_{1} \geq 0\}$, this is theorem 1.3.3. in \cite{Ho2}.}%
. We extend $A_{1}$ and $A_{2}$ as $\mathcal{C}^{2}$ compactly supported%
\footnote{Note that $A_{1}$ and $A_{2}$ do not necessarily agree on $\d \Omega$.}
functions in $\R^{n}$. 

We consider two geometrical optics solutions 
    $$ u_{j}(x,h)=e^{\frac{1}{h}(\phi_{j}+i\psi_{j})}(e^{\Phi_{j}}+hr_{j}(x,h)) $$
of the equations $\mathcal{L}_{A_{1},\bar{q}_{1}}u_{1}=0$ and  $\mathcal{L}_{A_{2},q_{2}}u_{2}=0$ constructed in the 
former section with phases  
\begin{align}
       \phi_{2}(x)&=-\phi_{1}(x)=\phi(x)=\log|x-x_{0}|,  \\ \nonumber
       \psi_{2}(x)&=\psi_{1}(x)=\psi(x)=d_{S^{n-1}}\Big(\frac{x-x_{0}}{|x-x_{0}|}, \omega \Big),  
\end{align}
defined on a neighborhood $\tilde{\Omega}$ of $\Omega$ (and $\omega$ varies in $\Gamma_{0}$), 
and where $\Phi_{1}$ and $\Phi_{2}$ are solutions of the equations
\begin{align}
\label{recovery:Phi}
    (\nabla \phi-i\nabla \psi) \cdot \nabla \Phi_{1}+i(\nabla \phi-i\nabla \psi)\cdot A_{1} 
    +\frac{1}{2} \Delta (\phi-i \psi) &=0 \\ \nonumber
    (\nabla \phi+i\nabla \psi) \cdot \nabla \Phi_{2}+i(\nabla \phi+i\nabla \psi)\cdot A_{2} 
    +\frac{1}{2} \Delta (\phi+i \psi) &=0.
\end{align}
The remainders $r_{j}$ are bounded independently of $h$ in $H^{1}_{\rm scl}$. Note that it implies the following estimate on $u_{j}$
\begin{align}
\label{recovery:estSol}
     \|e^{-\frac{\phi_{j}}{h}}u_{j}\|_{H^{1}_{\rm scl}}=\mathcal{O}(1).
\end{align}

By $w$ we denote the solution to the equation
    $$ \mathcal{L}_{A_{1},q_{1}}w=0, \quad w|_{\d \Omega}=u_{2}|_{\d \Omega} $$
so that $\mathcal{N}_{A_{1},q_{1}}(u_{2}|_{\d \Omega})=(\d_{\nu}w)|_{\d \Omega}+iA_{1}\cdot \nu \, u_{2}|_{\d \Omega}$.  
The assumption (\ref{recovery:hyp}) means that
   $$ \d_{\nu}(w-u_{2})(x) =0, \quad \forall x \in F_{\eps} $$
(here we use the fact (\ref{recovery:normComp}) that the normal components of the magnetic potentials coincide on the boundary).
Besides, we have
\begin{align}
\label{recovery:equ}
     \mathcal{L}_{A_{1},q_{1}}(w&-u_{2}) = -\mathcal{L}_{A_{1},q_{1}}u_{2}=
     (\mathcal{L}_{A_{2},q_{2}}-\mathcal{L}_{A_{1},q_{1}})u_{2} \\  \nonumber
     &=(A_{2}-A_{1})\cdot Du_{2}+D \cdot (A_{2}-A_{1}) u_{2}\\ \nonumber &\quad+(A_{2}^{2}-A_{1}^{2}+q_{2}-q_{1})u_{2}
\end{align}
hence we deduce
\begin{align}
\label{recovery:csqSol}
     \big(\mathcal{L}&_{A_{1},q_{1}}(w-u_{2})|u_{1}\big) \\  \nonumber &=\big((A_{2}-A_{1})\cdot Du_{2}|u_{1}\big)
     + \big(u_{2}|(A_{2}-A_{1})  \cdot Du_{1}\big)
     \\ \nonumber &\quad+\frac{1}{i}\underbrace{\big((A_{2}-A_{1})\cdot \nu \, u_{2}|u_{1}\big)_{\d \Omega}}_{=0} 
     +\big((A_{2}^{2}-A_{1}^{2}+q_{2}-q_{1})u_{2}|u_{1}\big).
\end{align}
The magnetic Green's formula gives
 \begin{align}
 \label{recovery:csqGreen}
      (\mathcal{L}_{A_{1},q_{1}}(w&-u_{2})|u_{1})\\ \nonumber 
      &= \underbrace{(w-u_{2}|\mathcal{L}_{A_{1},\bar{q}_{1}}u_{1})}_{=0}
      -\big((\d_{\nu}+iA_{1}\cdot \nu)(w-u_{2})|u_{1})_{\d\Omega}\\ \nonumber 
      &=-(\d_{\nu}(w-u_{2})|u_{1})_{\d\Omega \backslash F_{\eps}}
 \end{align}
 and combining (\ref{recovery:csqSol}) and (\ref{recovery:csqGreen}), we finally obtain
 \begin{multline}
 \label{recovery:theFormula}
     \int_{\d \Omega \backslash F_{\eps}} \d_{\nu}(u_{2}-w) \, \bar{u}_{1} \, d\sigma(x)=
     \int_{\Omega} (A_{2}^{2}-A_{1}^{2}+q_{2}-q_{1})u_{2} \, \bar{u}_{1} \, dx \\
     + \int_{\Omega} (A_{2}-A_{1}) \cdot (Du_{2} \, \bar{u}_{1} + u_{2} \, \overline{Du_{1}}) \, dx. 
\end{multline}

With our choice of $\phi_{2}=\phi=\log|x-x_{0}|$, we have
    $$ F_{\eps} \supset  F(x_{0})=\d \Omega_{-} \quad \textrm{ thus } \quad \d\Omega \backslash F_{\eps} \subset \d \Omega_{+} $$
 and moreover $\d_{\nu}\phi > \eps$ on $\d \Omega \backslash F_{\eps}$ therefore the modulus of the left-hand side in 
 (\ref{recovery:theFormula}) is bounded by
 \begin{align*}
      \frac{1}{\sqrt{\eps}}  \|\sqrt{\d_{\nu}\phi} \,  e^{-\frac{\phi}{h}} \d_{\nu}(u_{2}-w) \|_{\d \Omega_{+}} \times 
      \underbrace{\|e^{\Phi_{1}}+h r_{1}\|_{\d \Omega_{+}}}_{\leq \|e^{\Phi_{1}}\|_{\d\Omega}+\|r_{1}\|_{H^{1}_{\rm scl}}}
 \end{align*}
 which, in virtue of the Carleman estimate (\ref{Carleman:CarlEst}), is bounded by a constant times
 \begin{align*}
      \frac{1}{\sqrt{\eps}}\big(\sqrt{h} \|e^{-\frac{\phi}{h}}\mathcal{L}_{A_{1},q_{1}}(u_{2}-w)\|
      +\underbrace{\|\sqrt{-\d_{\nu}\phi} \,  e^{-\frac{\phi}{h}} \d_{\nu}(u_{2}-w) \|_{\d \Omega_{-}}}
      _{=0 \textrm{ because of }(\ref{recovery:hyp})}\big).
\end{align*}
In view of (\ref{recovery:equ}) and of (\ref{recovery:estSol})  the former expression is 
$\mathcal{O}(h^{-\frac{1}{2}})$.  Therefore we can conclude that the right-hand side of (\ref{recovery:theFormula}) is 
$\mathcal{O}(h^{-\frac{1}{2}})$.  This constitutes an important difference with \cite{KSU}, where%
\footnote{See (5.16) and the subsequent lines in \cite{KSU}.}
the corresponding term was $\mathcal{O}(h)$.  

More directly, the first and second right-hand side terms of (\ref{recovery:theFormula}) are respectively $\mathcal{O}(1)$ and
$\mathcal{O}(h^{-1})$ as may be seen from (\ref{recovery:estSol}). It turns out that the information obtained when disregarding 
the bounded term is enough to recover the magnetic field. We multiply (\ref{recovery:theFormula}) by $h$ and let $h$ tend to $0$:    
\begin{align*}
    \lim_{h \rightarrow 0} \int_{\Omega} \big((A_{2}-A_{1})\cdot h Du_{2} \, \bar{u}_{1}
     + u_{2} \, (A_{2}-A_{1})  \cdot h \overline{Du_{1}} \big) \, dx =0. 
\end{align*}   
Using the explicit form of the solutions $u_{1}$ and $u_{2}$, this further means
\begin{align}
\label{recovery:limit}
   \int_{\Omega} (A_{2}-A_{1}) \cdot (\nabla \phi+i\nabla \psi) e^{\bar{\Phi}_{1}+\Phi_{2}} \, dx =0. 
\end{align} 
Adding the complex conjugate of the first line of (\ref{recovery:Phi}) to the second line, we see that
\begin{align*}
    (\nabla \phi+i\nabla \psi) \cdot \big(\nabla \bar{\Phi}_{1}+\nabla \Phi_{2}+ i(A_{2}-A_{1})\big)
    +\Delta (\phi+i \psi) =0
\end{align*}
this implies that
\begin{align}
\label{recovery:PhiEq}
   (D+ A_{2}-A_{1})  \cdot (\nabla \phi+i\nabla \psi) (e^{\bar{\Phi}_{1}+\Phi_{2}})=0.
\end{align}

As observed in remark \ref{WKB:remhol},  in the expression for $u_{2}$, we may replace
$e^{\Phi_{2}}$ by $e^{\Phi_{2}}g$ if $g$ is a solution of $(\nabla \phi + i \nabla \psi) \cdot \nabla g =0$. 
Then (\ref{recovery:limit}) can be replaced by
\begin{align*}
   \int_{\Omega} (A_{2}-A_{1}) \cdot (\nabla \phi+i\nabla \psi) e^{\bar{\Phi}_{1}+\Phi_{2}}g \, dx =0. 
\end{align*} 
From equation (\ref{recovery:PhiEq}), we see that we can replace $A_{2}-A_{1}$ by $i \nabla$ in the former equality
\begin{align}
   \int_{\Omega} g(x) \, \nabla \cdot \big(e^{\bar{\Phi}_{1}+\Phi_{2}} (\nabla \phi+i\nabla \psi)\big)  \, dx =0
\end{align} 
for all functions $g$ such that $(\nabla \phi + i \nabla \psi) \cdot \nabla g =0$ on $\Omega$.
\begin{rem}
    An integration by parts gives
    \begin{align*}
         \int_{\d \Omega} (\d_{\nu}\phi+i\d_{\nu}\psi) e^{\bar{\Phi}_{1}+\Phi_{2}}g \, d\sigma 
         -\int_{\Omega} e^{\bar{\Phi}_{1}+\Phi_{2}} \underbrace{(\nabla \phi+i\nabla \psi) \cdot \nabla g}_{=0} \, dx =0
    \end{align*} 
    hence we have $\displaystyle{\int_{\d \Omega} (\d_{\nu}\phi+i\d_{\nu}\psi) e^{\bar{\Phi}_{1}+\Phi_{2}}g \, d\sigma =0}$.
 \end{rem}
\end{section}
%
%%%%%%%%%%%%%%%%%%%%%%%%%%%%%%%%%%%%%%%%%%%%%%%%%%%%%%%%%%%%%%
%%%%%%%%%%%%%%%%%%%%%%%%%%%%%%%%%%%%%%%%%%%%%%%%%%%%%%%%%%%%%%
%%                                           		      	 FINAL ARGUMENTS                                                                                         %%
%%%%%%%%%%%%%%%%%%%%%%%%%%%%%%%%%%%%%%%%%%%%%%%%%%%%%%%%%%%%%%
%%%%%%%%%%%%%%%%%%%%%%%%%%%%%%%%%%%%%%%%%%%%%%%%%%%%%%%%%%%%%%
%
\begin{section}{Moving to the complex plane}
\label{complex}
In this section, we follow remark \ref{WKB:remzvar} and choose to work in the cylindrical coordinates. 
Let us be more precise: suppose that $x_{0}=0 \notin \ch(\Omega)$ and that we have picked $\omega \in S^{n-1} \backslash (\Gamma  
\cup \check{\Gamma})$ with the notations of section \ref{WKB}. After a rotation, we assume that $\omega=(1,0,\dots,0)$, therefore
we have
    $$ \Omega \Subset \tilde{\Omega}  \subset \{x \in \R^{n}: x' \neq 0\}. $$
We choose the following cylindrical coordinates
    $$ t=x_{1}, \quad r=|x'|>0, \quad \theta = \Theta(x)=\frac{x'}{|x'|} \in S^{n-2}. $$
By Sard's theorem, the set of critical values of $\Theta : \Omega \mapsto S^{n-2}$ is of measure $0$,
therefore the set $\Omega_{\theta_{0}}=\Theta^{-1}(\theta_{0})=\{x \in \Omega : x'=r \theta_{0}, r>0\}$
is an open set with smooth boundary for almost every $\theta_{0}$ in $\Theta(\Omega)$.
The result obtained in the former section reads
\begin{align*}
   \int_{S^{n-2}} \iint_{\Omega_{\theta}} g(x) \, \nabla_{x} \cdot \big(e^{\bar{\Phi}_{1}+\Phi_{2}} (\nabla_{x} \phi+i\nabla_{x} \psi)
   \big) \, r^{n-2}dr  dt \, d\theta =0. 
\end{align*} 
and taking $g=g_{1}(t,r) \otimes g_{2}(\theta)$ and varying $g_{2}$ leads to 
\begin{align}
\label{complex:main}
    \iint_{\Omega_{\theta}} g(t,r) \, \nabla_{x} \cdot \big(e^{\bar{\Phi}_{1}+\Phi_{2}} (\nabla_{x} \phi+i\nabla_{x} \psi)\big) 
    \, r^{n-2}dr  dt  =0
\end{align} 
for any function $g$ such that $(\nabla_{x}\phi+i\nabla_{x}\psi) \cdot \nabla_{x}g=0$ on $\Omega_{\theta}$, and this for almost every $\theta$.

Now we consider the complex variable  $z=t+ir \in \C_{+}=\{w \in \C: \im w>0\}$ and write our results in this setting.
Let us recall the results of the computations made in remarks \ref{WKB:remzvar} and \ref{WKB:remzvarbis}:
 \begin{align}
         \nabla_{x}(\phi + i\psi) \cdot \nabla_{x} =\frac{2}{z}\frac{\d}{\d \bar z}
          \quad &\textrm{and} \quad \Delta_{x} (\phi+i\psi)=\frac{(n-2)i}{rz} \\ \nonumber
          \textrm{ and thus } \quad  \nabla_{x} \cdot \nabla_{x} \circ (\phi + i\psi) &=\frac{2}{z}\frac{\d}{\d \bar z}+\frac{(n-2)i}{r z}. 
\end{align}
Similarly, the functions $\Phi_{j}$ satisfy
\begin{align}
\label{complex:PhiEq}
      \frac{\d }{\d \bar{z}}(\bar{\Phi}_{1}+\Phi_{2})=\frac{(n-2)}{(z-\bar{z})} + \frac{1}{2} (A_{1}-A_{2}) \cdot (e_{1}+ie_{r})
\end{align}
Finally, (\ref{complex:main}) reads
\begin{align}
\label{complex:mainComplex}
   \iint_{\Omega_{\theta}} g(z) \frac{1}{z}\Big(\frac{\d}{\d \bar z}-\frac{(n-2)}{(z-\bar{z})}\Big) (e^{\bar{\Phi}_{1}+\Phi_{2}})
   (z-\bar{z})^{n-2} \, d\bar{z} \wedge dz=0,
\end{align} 
for any  $g \in \mathcal{H}(\Omega_{\theta})$.  Replacing the holomorphic function $g/z$ on $\Omega_{\theta}$ by $g$, we can drop
the factor $1/z$.

If $g$ is a holomorphic function, we have%
\footnote{The result of this computation is the transcription of the fact that the formal adjoint of $(\nabla \phi+i\nabla \psi)\cdot \nabla$
is $\nabla \cdot (\nabla \phi+i\nabla \psi)$ in the complex setting, where  the measure is $r^{n-2} d\bar{z} \wedge dz$.}
\begin{align*}
    d\big(e^{\bar{\Phi}_{1}+\Phi_{2}} g(z)& (z-\bar{z})^{n-2} dz \big)  \\
    &=\frac{\d}{\d \bar{z}} \big( (z-\bar{z})^{n-2}  e^{\bar{\Phi}_{1}+\Phi_{2}}\big)  g(z)  d\bar{z} \wedge dz  \\
    &= \big(\frac{\d}{\d \bar{z}}-\frac{n-2}{z-\bar{z}}\big) \big( e^{\bar{\Phi}_{1}+\Phi_{2}}\big)  g(z)    
    (z-\bar{z})^{n-2} d\bar{z} \wedge dz 
\end{align*}
therefore the Stokes' formula implies
\begin{multline*}
   \iint_{\Omega_{\theta}} g(z) \Big(\frac{\d}{\d \bar z}-\frac{(n-2)}{(z-\bar{z})}\Big) (e^{\bar{\Phi}_{1}+\Phi_{2}})
   (z-\bar{z})^{n-2} \, d\bar{z} \wedge dz \\ = \int_{\d \Omega_{\theta}} g(z) e^{\bar{\Phi}_{1}+\Phi_{2}}
   (z-\bar{z})^{n-2} \, dz.
\end{multline*}
Together with  (\ref{complex:mainComplex}) this gives
\begin{align}
\label{complex:boundary}
    \int_{\d \Omega_{\theta}} g(z) e^{\bar{\Phi}_{1}+\Phi_{2}}(z-\bar{z})^{n-2} \, dz =0
\end{align}
for any  $g \in \mathcal{H}(\Omega_{\theta})$. 

\begin{lem}
\label{complex:HolTrace}
    There exists a non-vanishing holomorphic function $F$ on $\Omega_{\theta}$, continuous on $\bar{\Omega}_{\theta}$, 
    whose restriction to $\d\Omega_{\theta}$ is equal to $(z-\bar{z})^{n-2}e^{\bar{\Phi}_{1}+\Phi_{2}}$.
\end{lem}
\begin{proof}
    We denote $f(z)=(z-\bar{z})^{n-2}e^{\bar{\Phi}_{1}+\Phi_{2}}$ and consider the Cauchy integral operator
           $$ C(f)(z)=\frac{1}{2\pi i} \int_{\d \Omega_{\theta}} \frac{f(\zeta)}{\zeta-z} \, d\zeta, \quad \forall z \in \C \backslash \d \Omega_{\theta}. $$
    The function $C(f)$ is holomorphic inside and outside $\Omega_{\theta}$ and the Plemelj-Sokhotski-Privalov formula reads on the boundary
    \begin{align}
    \label{complex:plemelj}
          \lim_{\substack{z \rightarrow z_{0} \\ z \in \Omega_{\theta}}}C(f)(z)-\lim_{\substack{z \rightarrow z_{0} \\ z \notin \Omega_{\theta}}}C(f)(z)
           = f(z_{0}), \quad \forall z_{0} \in \d \Omega_{\theta}.
    \end{align}
    The function $\zeta \rightarrow (\zeta-z)^{-1}$ is holomorphic on $\Omega_{\theta}$ when $z \notin \Omega_{\theta}$ hence
    (\ref{complex:boundary}) implies that $C(f)(z)=0$ when $z \notin \Omega_{\theta}$. The second limit in (\ref{complex:plemelj})
    is $0$, thus $F=C(f)$ is a holomorphic function on $\Omega_{\theta}$ whose restriction to the boundary agrees with $f$.
     
    It remains to prove that $F$ does not vanish on $\Omega_{\theta}$.  This is clear by the argument principle since
         $$ \vararg_{\d \Omega_{\theta}}F =\vararg_{\d \Omega_{\theta}}f  = 0 $$
    and $F$ is holomorphic.         
\end{proof}
In particular, with the former function, we have on the boundary
    $$ \frac{1}{F(z)}(z-\bar{z})^{n-2} e^{\bar{\Phi}_{1}+\Phi_{2}} =1, \quad \forall z \in \d\Omega_{\theta} $$
which implies
   $$ \bar{\Phi}_{1}+\Phi_{2} + \log (z-\bar{z})^{n-2} = \log F(z), \quad \forall z \in \d\Omega_{\theta} $$
with $\log F$ a holomorphic function on $\Omega_{\theta}$, and therefore
\begin{align*}
    \int_{\d \Omega_{\theta}} g(z) \Big(\bar{\Phi}_{1}+\Phi_{2}+ \log (z-\bar{z})^{n-2}\Big) \, dz=0.
\end{align*}
An application of  Stokes' formula gives
\begin{align*}
    \iint_{\Omega_{\theta}} g(z) \Big(\frac{\d}{\d \bar{z}} (\bar{\Phi}_{1}+\Phi_{2})-
    \frac{n-2}{z-\bar{z}} \Big) \, d\bar{z} \wedge dz=0
\end{align*}
hence using the equation (\ref{complex:PhiEq}),  this implies
\begin{align}
\label{complex:final}
    \iint_{\Omega_{\theta}} g(z) (A_{1}-A_{2}) \cdot (e_{1}+ie_{r})   \,   d\bar{z} \wedge dz  =0.
\end{align}
With $g=1$, the former equality reads
\begin{align*}
    \iint_{\Omega_{\theta}} (t+ir) (\nabla_{x}\phi+i\nabla_{x} \psi) \cdot (A_{1}-A_{2}) \, dt \, dr =0.
\end{align*}

Denote by $P_{\theta}=\mathop{\rm span}(\omega,e_{r})$ the plane along the axis directed by $\omega=(1,0,\dots,0)$,
and by $P^{+}_{\theta}$ the half plane where $x \cdot e_{r} >0$, then
$\Omega_{\theta}=\Omega \cap \{x=(x_{1},x') \in \R^{n}: x'=r \theta, r>0 \}=\Omega \cap P^{+}_{\theta}$. 
Let $\pi_{\theta}$ be the projection on $P_{\theta}$ and $d\lambda_{\theta}$ the measure on the plane, then (\ref{complex:final})
(with $g=1$) implies
    $$ \int_{P_{\theta}\cap \Omega} \pi_{\theta} (A_{1}-A_{2})\, d\lambda_{\theta} =0, $$
for almost every $\theta \in S^{n-2}$, hence for all $\theta \in S^{n-2}$ by continuity. 
The former may be rephrased under the form
    $$ \int_{x_{0}+P}  \xi \cdot \big(1_{\Omega}(A_{1}-A_{2})\big)\, d\lambda_{P} =0, \quad \forall \xi \in P $$
for all linear planes $P$ containing $\omega=(1,0,\dots,0)$.  We can also let $x_{0}$ vary in a small neighborhood of $0 \notin \ch(\Omega)$ and $\omega$ vary in a neighborhood $\Gamma_{0}$ of $(1,0,\dots,0)$ on the sphere $S^{n-1}$.

\begin{lem}
\label{complex:helgason}
    Let $A$ be a $\mathcal{C}^{1}$ vector field on $\bar{\Omega}$. If
    \begin{align}
    \label{complex:radon}
          \int_{P \cap \Omega}  \xi \cdot A \, d\lambda_{P} =0, \quad \forall \xi \in T_{x}(P)
    \end{align}
    for all planes $P$ such that $d((0,e_{1}),T(P))<\delta$ then $d A=0$ on $\Omega$.
\end{lem}
The proof of this lemma is based on the following microlocal version of Helgason's support theorem.
\begin{thm}
\label{complex:microHelgason}
    Let $f \in \mathcal{C}^{0}(\R^{n})$, suppose that the Radon transform of $f$ satisfies
          $$ \mathcal{R}f(H) =\int_{H} f \, d\lambda_{H} = 0 $$
    for all hyperplanes $H$ in some neighbourhood of a hyperplane $H_{0}$ then
          $$ N^{*}(H_{0}) \cap {\rm WF}_{a}(f) = \varnothing $$
    where $N^{*}(H_{0})$ cenotes the conormal bundle of $H_{0}$.
\end{thm}
The proof of this result may be found in \cite{Bo} (see proposition 1) or in \cite{Ho} (see section 6).
We will also need the microlocal version of Holmgren's theorem (see \cite{Ho}, section 1 or \cite{Ho2} section 8.5).
\begin{thm}
\label{complex:microHolmgren}
    Let $f \in \mathcal{E}'(\R^{n})$ then we have 
        $$  N(\supp f) \subset {\rm WF}_{a}(f) $$
    where $N(\supp f)$ is the normal set of the support of $f$.
\end{thm}
These two results may be combined to provide a proof of Helgason's support theorem (see \cite{Bo} and \cite{Ho}).
We also refer to the book \cite{He} for a review on Radon transforms. 
\begin{proof}[Proof of lemma \ref{complex:helgason}]
     Let $\chi \in \mathcal{C}^{\infty}_{0}(|x|<\frac{1}{2})$ and $\chi_{\eps}=\eps^{-n}\chi(\cdot/\eps)$ be a standard regularization,  one has   
                    $$ \int_{P} \xi \cdot (\chi_{\eps}*1_{\Omega}A) \, d\lambda_{P} = \eps^{-n} \int \chi\big(\frac{y}{\eps}\big) 
                         \Big(\int_{(-y+P) \cap \Omega} \xi \cdot A \, d\lambda_{-y+P}\Big) \, dy = 0 $$
     when $d((0,e_{1}),T(P))<\delta-\eps$. Therefore it suffices to prove the result when $\Omega=\R^{n}$ and $A \in \mathcal{C}^{\infty}_{0}(\R^{n};\R^{n})$
     since $d(\chi_{\eps}*1_{\Omega}A)$ tends to $dA$ as a distribution when $\eps$ tends to $0$.
    
    Our first step is to prove that
    \begin{align}
    \label{complex:restr}
         \iota_{H}^{*}dA=0
    \end{align}
    for any subspace $H \subset \R^{n}$ of dimension 3 such that $d((0,e_{1}),T(H))<\delta$. Here $\iota_{H}$ denotes the injection of $H$ in $\R^{n}$.
    For any plane $P \subset H$ such that $d((0,e_{1}),T(P))<\delta$, we have
          $$ \int_{P} \langle dA , \xi \wedge \eta \rangle \, d\lambda_{P} = 
               \frac{d}{dt}\Big( \int_{t\xi+P} \eta \cdot A \, d\lambda_{P} -  \int_{t\eta+P} \xi \cdot A \, d\lambda_{P} \Big)_{t=0} $$
    when $\xi,\eta \in T_{x}(H)$.  The space $H$ is of dimension 3 so we can assume that either $\eta$ or $\xi$ belongs to $T_{x}(P)$, thus 
    the former expression is zero because of (\ref{complex:radon}) and of the fact that whenever $\eta \in T_{x}(P)$
         $$  \int_{t\eta+P} \xi \cdot A \, d\lambda_{P} \quad \textrm{ is constant.} $$
         
    Therefore, if $\mathcal{R}_{H}$ denotes the Radon transform in $H$, we obtain
         $$ \mathcal{R}_{H}\big(\langle \iota_{H}^{*}dA,\xi \wedge \eta \rangle\big)(P) = 0 $$
    for any plane $P\subset H$ such that $d((0,e_{1}),T(P))<\delta$ and for any $\xi,\eta \in T_{x}(H)$.
    Combining theorems \ref{complex:microHelgason} and \ref{complex:microHolmgren} we obtain
         $$ N^{*}(P) \cap N(\supp \langle \iota_{H}^{*}dA,\xi \wedge \eta \rangle) = \varnothing $$
    for any plane $P\subset H$ such that $d((0,e_{1}),T(P))<\delta$. This gives (\ref{complex:restr})
    since such a family of planes sweep other $\Omega \cap H$ and the support of $A$ is on one side of at least one such plane.
    
    The result (\ref{complex:restr}) implies in particular that
          $$ \langle dA(x), \xi \wedge \eta \rangle = 0, \quad \forall x \in \R^{n}, \forall (\xi,\eta) \in S^{n}\times \R^{n}, \, |\xi-e_{1}|<\delta $$   
      and therefore $dA=0$  by linearity. 
\end{proof}
\end{section}
%
%%%%%%%%%%%%%%%%%%%%%%%%%%%%%%%%%%%%%%%%%%%%%%%%%%%%%%%%%%%%%%
%%%%%%%%%%%%%%%%%%%%%%%%%%%%%%%%%%%%%%%%%%%%%%%%%%%%%%%%%%%%%%
%%                                               		          	 POTENTIAL                                                                                             %%
%%%%%%%%%%%%%%%%%%%%%%%%%%%%%%%%%%%%%%%%%%%%%%%%%%%%%%%%%%%%%%
%%%%%%%%%%%%%%%%%%%%%%%%%%%%%%%%%%%%%%%%%%%%%%%%%%%%%%%%%%%%%%
%
\begin{section}{Recovering the potential}
\begin{proof}[End of the proof of theorem \ref{intro:mainThm}]
    Applying this lemma, we finally obtain
         $$ dA_{1}=dA_{2} \quad \textrm{ on } \Omega $$
    therefore the difference of the two potentials is a gradient $A_{1}-A_{2}=\nabla \Psi$ (recall that $\Omega$ is simply connected).
    The identity (\ref{complex:final}) now reads
    \begin{align*}
         \iint_{\Omega_{\theta}} g(z) \d_{\bar{z}}\Psi(z,\theta)  \,   d\bar{z} \wedge dz  =0
    \end{align*}
    for any holomorphic function $g$ on $\Omega_{\theta}$ and by Stokes' theorem we get
    \begin{align*}
         \int_{\d\Omega_{\theta}} g(z) \Psi(z,\theta)  \,   dz  =0.
    \end{align*}
    Reasoning as in the beginning of lemma \ref{complex:HolTrace}, there exists a holomorphic function $\tilde{\Psi} \in 
    \mathcal{H}(\Omega_{\theta})$ such that $\tilde{\Psi}|_{\d\Omega_{\theta}}=\Psi|_{\d\Omega_{\theta}}$. 
    Now $\Psi$ is real-valued, and since $\tilde{\Psi}$ is real-valued on $\d\Omega_{\theta}$ and harmonic, it is real-valued everywhere. 
    The only real-valued holomorphic functions are the constant ones, so $\tilde{\Psi}$ and hence $\Psi$ is constant on $\d\Omega_{\theta}$. 
    Varying $\theta$ and also slightly $x_{0}$ and $\omega$, we get that  $\Psi$  is constant on the boundary $\d\Omega$.
    We may assume then that $\Psi|_{\d\Omega}=0$.
    
    By a gauge transformation, we may assume that $\Psi=0$, thus $A_{1}=A_{2}$.  We could almost directly apply the result in \cite{KSU} 
    to recover the identity of the two potentials $q_{1}=q_{2}$,  if it were not for the presence of the two magnetic potentials in the equations.
    Instead we go back to the limit induced by (\ref{recovery:theFormula}). The second right-hand side term is now zero. The
    left-hand side is now $\mathcal{O}(\sqrt{h})$ since the $\mathcal{O}(h^{-1})$ term in (\ref{recovery:equ}) is zero  and we
    can reproduce the arguments given after (\ref{recovery:theFormula}). Therefore we obtain 
    \begin{align}
    \label{potential:theLim}
         \lim_{h \rightarrow 0} \int_{\Omega} (q_{2}-q_{1})u_{2} \bar{u}_{1}\, dx =0 
    \end{align}
    thus
         $$ \int_{\Omega} (q_{2}-q_{1}) e^{\bar{\Phi}_{1}+\Phi_{2}} \, dx = 0. $$
    As observed in remark \ref{WKB:remhol},  we may replace $e^{\Phi_{2}}$ by $e^{\Phi_{2}}g$ 
    if $g$ is a solution of $(\nabla \phi + i \nabla \psi) \cdot \nabla g =0$. Then the former can be replaced by
         $$ \int_{\Omega} (q_{2}-q_{1}) e^{\bar{\Phi}_{1}+\Phi_{2}} g(x) \, dx = 0. $$
    Moving to the complex plane, as in section \ref{complex}, this reads
         $$ \iint_{\Omega_{\theta}} (q_{2}-q_{1})g(z) e^{\bar{\Phi}_{1}+\Phi_{2}} (z-\bar{z})^{n-2} \, d\bar{z} \wedge dz = 0 $$
    for any holomorphic function $g$ on $\Omega_{\theta}$.  But the transport equation (\ref{complex:PhiEq}) now reads 
         $$ \frac{\d}{\d \bar{z}} \big((z-\bar{z})^{n-2}e^{\bar{\Phi}_{1}+\Phi_{2}}\big)=0 $$
    therefore if we take $g=(z-\bar{z})^{-n+2}e^{-\bar{\Phi}_{1}-\Phi_{2}}$, we obtain
         $$ \iint_{\Omega_{\theta}} q(t,r,\theta) \, dt \, dr = 0 $$
    with $q=q_{1}-q_{2}$. As in section \ref{complex}, varying $x_{0}$ and $\omega$, this may be interpreted as
         $$ \int_{P} 1_{\Omega}q \, d\lambda_{P} = 0 $$
    for any plane $P$ such that $d((0,e_{1}),T(P))<\delta$. This implies that for any subspace $H$ of dimension 3 
    such that $d((0,e_{1}),T(H))<\delta$ we have
        $$ \mathcal{R}_{H}(1_{\Omega}q)(P) = \int_{P} 1_{\Omega}q \, d\lambda_{P} = 0 $$
    for any plane such that $d((0,e_{1}),T(P))<\delta$.  Regularizing and applying theorems \ref{complex:microHelgason} and \ref{complex:microHolmgren}
    as in the proof of lemma \ref{complex:helgason}, we get
        $$ N^{*}(P) \cap N(\supp (1_{\Omega}q)|_{H}) = \varnothing $$
    and therefore $1_{\Omega}q=0$ on $H$ leading to $q=0$ on $\Omega$.
    This ends the proof of theorem \ref{intro:mainThm}.  
 \end{proof}
\end{section}
\begin{acknow}
   The work of Carlos Kenig and Gunther Uhlmann was partially supported by NSF. Johannes Sj\"ostrand wishes to acknowledge the 
hospitality of the Department of Mathematics of the University of Washington, and David Dos Santos Ferreira  the 
hospitality of the Department of Mathematics of the University of Chicago.
\end{acknow}
%
%%%%%%%%%%%%%%%%%%%%%%%%%%%%%%%%%%%%%%%%%%%%%%%%%%%%%%%%%%%%%%
%%%%%%%%%%%%%%%%%%%%%%%%%%%%%%%%%%%%%%%%%%%%%%%%%%%%%%%%%%%%%%
%%                                                                         BIBLIOGRAPHY                                                                                          %%
%%%%%%%%%%%%%%%%%%%%%%%%%%%%%%%%%%%%%%%%%%%%%%%%%%%%%%%%%%%%%%
%%%%%%%%%%%%%%%%%%%%%%%%%%%%%%%%%%%%%%%%%%%%%%%%%%%%%%%%%%%%%%
%

%
%
\end{document}